# RIEMANN HYPOTHESIS: ARCHITECTURE OF A CONJECTURE 'ALONG' THE LINES OF PÓLYA. FROM TRIVIAL ZEROS AND HARMONIC OSCILLATOR TO INFORMATION ABOUT NON-TRIVIAL ZEROS OF THE RIEMANN ZETA-FUNCTION

STEFANO BELTRAMINELLI, DANILO MERLINI, AND SERGEY SEKATSKII

**Abstract.** We propose an architecture of a conjecture concerning the Riemann Hypothesis in the form of an 'alternative' to the Pólya strategy: we construct a Hamiltonian $\hat{H}_{Polya}$ whose spectrum coincides exactly with that of the Harmonic Oscillator Hamiltonian $\hat{H}_{osc}$ if and only if the Riemann Hypothesis holds true. In other words, it can be said that we formulate the Riemann Hypothesis by means of a non-commutative structure
on the real axis, viz., that of the Harmonic Oscillator, by an equation of the type $\hat{H}_{Polya}(\hat{H}_{osc}) = \hat{H}_{osc}$: the Harmonic Oscillator operator, if viewed as an argument of $\hat{H}_{Polya}$, reproduces itself.

## 1. Introduction

Many years ago Landau asked to Pólya if he knows a physical reason why the Riemann Hypothesis (RH; see e.g. [1] for definition and discussion of the general properties of Riemann $\varsigma$-function) should be true [2]. The 'answer' of Pólya was that RH would be true if someone is able to find a 'Hamiltonian' (a self-adjoint unbounded operator $\hat{T}$) whose eigenvalues (necessarily real, of course) coincide with, say, $i(1/2 - \rho_k)$, where $\rho_k$ are the non-trivial Riemann $\varsigma$-function zeros.

Many researchers constructed systems of similar nature by mathematical and physical means, and gave interesting and important contributions to the problem. As an example, Berry, Keating and others proposed for $\hat{T}$ a one-dimensional operator of the form $\hat{T} = x \cdot p$, to be quantized in some special space. See, among others, the results given in this direction in [3-8].

In this work we attempt to formulate the problem, still in one dimension, by means of an operator $\hat{H}_{osc}$, i.e. that of the Harmonic Oscillator. Instead of the search for an operator $\hat{T}$ possibly having $i(1/2 - \rho_k)$ values as eigenvalues, here we start from an operator $\hat{H}_{osc}$ connected with the trivial zeros of the Riemann $\varsigma$-function. These limited from above and regularly sequenced trivial zeros clearly can be viewed as the eigenvalues of some Harmonic Oscillator operator, and our main idea can be formulated in the following way. Starting from an appropriately defined $\hat{H}_{osc}$ operator, whose spectrum describes all the trivial zeros, we are trying to modify it and to 'enlarge' its properties in such a way that, besides the trivial zeros, this new $\hat{H}_{Polya}$ operator should also 'catch' some useful information about non-trivial zeros. This new operator is constructed starting from an appropriate integral equality for the function $\varsigma'(z)/\varsigma(z)$ (see below) which is closely connected with a number of integral equalities involving the function $\ln(\varsigma(z))$ and shown by us earlier to be equivalent to the Riemann Hypothesis [9-11]. An appropriate Taylor expansion of analytic functions having the operator $\hat{H}_{osc}$ as an argument, similarly to the well-known construction of an operator $\exp(-a\hat{H})$ starting from an operator $\hat{H}$ treated in many textbooks, is exploited. On RH, the so obtained operator $\hat{H}_{Polya}$ has exactly the same spectrum as the initial $\hat{H}_{osc}$ operator, but it has necessary different spectrum if RH fails, and in this sense we are speaking here about an 'architecture along the lines of Pólya'.

Thus the Hilbert-Pólya conjecture gives ideas about a possible 'similar' conjecture and a possible reformulation of the RH by means of the trivial zeros instead of the non-trivial zeros (which apparently are more difficult to treat). One may "suspect" that $\hat{H}_{Polya}$ is more artificial then the hypothetical $\hat{T}$, but it should be remembered that the Harmonic Oscillator (quantum or classical) is a prototype of fundamental physical systems, and that integers are 'already quantized'. We do not invoke more about quantum mechanics but the connection of our approach with quantum mechanics seems to be quite strong and concrete.

Briefly, in our formulation, the Hamiltonian of a Harmonic Oscillator is the solution of an implicit equation involving the $\varsigma$-function where the Harmonic Oscillator itself is the argument of this Zeta function. The trivial zeros, i.e. the Harmonic Oscillator spectrum, thus lies in the very core of our 'alternative' to the Hilbert-Pólya conjecture formulated via a simple physical system. The paper is primarily of heuristic nature and we do not pretend on rigorous consideration of all related questions.

The first version of this work has been presented as a poster at the Quantum Mechanics, Operator Theory and the Riemann Zeta function Conference in Benasque, Spain, 17 – 23 June 2012.

## 2. An operator illustrating a modified Hilbert – Pólya approach

Following our previous papers [8 - 10], let us consider a rectangular contour $C$ with the vertices $b-i$, $b+iX$, $X+iX$, $X-i$, introduce a function $g(z) = i\exp(ia(z-b))$ and consider the contour integral $\int_C \frac{\varsigma'}{\varsigma}(z)g(z)dz$ (actually, instead of the line *(-i, X-i)* any other line *(-ci, X-ci)* can be used; *b* is real not equal to 1 and -2, -4, -6…; *a, c, X* are real positive). Along the line *(b+iX, X+iX)* we have $g(z) = O(e^{-aX})$ while along the line *(X-i, X+iX)* $\varsigma'(z)/\varsigma(z) = O(2^{-X})$, hence in the limit $X \to \infty$ the contour integral is equal to $\int_{-1}^{\infty} \frac{\varsigma'}{\varsigma}(b+it)e^{-at}dt + i\int_{0}^{\infty} \frac{\varsigma'}{\varsigma}(x+b-i)\exp(i(a(x-i)))dx$ and, due to Cauchy residue theorem, also is equal to the sum over poles of $\varsigma'(z)/\varsigma(z)$ function:

$$\int_C \frac{\varsigma'}{\varsigma}(z)i\exp(ia(z-b))dz = 2\pi i \sum_{\rho, \sigma_k > b} l_k e^{-at_k}(i\cos(a(\sigma_k - b)) - \sin(a(\sigma_k - b))) +$$

$$2\pi i \sum_{k=1}^{k<[-b/2]}(i\cos(a(-2k-b)) - \sin(a(-2k-b))) - 2\pi i(i\cos(a(1-b)) - \sin(a(1-b)))$$

In this sum one easily recognizes first the sum over non-trivial and then over trivial zeros, and finally the contribution of the simple pole of the Riemann $\varsigma$-function lying at *z=1*. Let us take $a = \pi$ and the sequence $b_n = -2n + 1/4$, here *n* is a positive integer (other choices of *a, $b_n$* are possible).

With such a choice we have for any *n*:

$$2\pi i \sum_{k=1}^{k<[-b_n/2]} (i\cos(\pi(-2k-b)) - \sin(\pi(-2k-b))) = \pi\sqrt{2}(n-1)i(i+1) = \pi\sqrt{2}(n-1)(-1+i).$$

For the sum over non-trivial zeros, provided all zeros lye on the critical line, we have $\pi\sqrt{2}i\sum_{\rho} l_k e^{-\pi t_k}(i+1) = \pi\sqrt{2}\sum_{\rho} l_k e^{-\pi t_k}(-1+i)$. It is easy to see that if there are paired zeros with $\sigma = 1/2 \pm \chi$ with positive $\chi \neq 0$, instead of the term $\pi\sqrt{2}l_k e^{-\pi t_k}(-1+i)$, in the sum over non-trivial zeros we have $\pi 2\sqrt{2}l_k e^{-\pi t_k}(-1+i)\cos(\pi\chi)$ so, due to $\pi\chi < \pi/2$ (there are no non-trivial zeros outside the critical strip), for both imaginary and real part the module of such a contribution is smaller than it would be for the case when the RH hypothesis holds. In other worlds, the sum $\Sigma_\pi = \sum_{\rho} l_k e^{-\pi t_k} \cos(\pi(1/2 - \sigma_k))$ is a *unique* characteristic of the Riemann $\varsigma$-function in a sense that the abscissa of no one zero can be changed without the change of the value of $\Sigma_\pi$. Or let us put it differently: provided we somehow compute the values of $\Sigma_\pi$ and $\Sigma_{\pi,RH} = \sum_{\rho} l_k e^{-\pi t_k}$, we then may use them to test the RH: RH holds if and only if $\Sigma_\pi = \Sigma_{\pi,RH}$. Unfortunately, at least at the current state of research, we are unable to efficiently calculate both these values. However, we believe that the very existence of such a criterion makes the subsequent consideration sufficiently interesting.

Contribution of the pole at *z=1* for our choice of *a* and $b_n$ is $\sqrt{2}\pi(-1-i)$. We thus have that

$$\int_{-1}^{\infty} \frac{\varsigma'}{\varsigma}(-2n+1/4+it)e^{-\pi t}dt + i\int_{0}^{\infty} \frac{\varsigma'}{\varsigma}(x-2n+1/4-i)\exp(i(\pi(x-i)))dx =$$
$$\pi\sqrt{2}(n-1)(-1+i) + \sqrt{2}\pi\Sigma_\pi(-1+i) + \sqrt{2}\pi(-1-i) \quad (1)$$

which we rewrite as

$$\frac{-1-i}{2\sqrt{2\pi}}\int_{-1}^{\infty} \frac{\varsigma'}{\varsigma}(-2(n+1/2)+5/4+it)e^{-\pi t}dt +$$
$$\frac{(-1-i)e^\pi}{2\sqrt{2\pi}}\int_{0}^{\infty} \frac{\varsigma'}{\varsigma}(y-2(n+1/2)+5/4)(i\cos(\pi y)-\sin(\pi y))dy = n+1/2+(\Sigma_\pi - 3/2 + i) \quad (2).$$

Here one can easily see that this rigorously proven equality contains *n+1/2* both as an "argument" of the functions involved and as an "answer". Moreover, no other "argument" of this equation apart from *n+1/2* where $n=0, 1, 2\ldots$ satisfies it (this is easy to see because $p=2$ is a period of the function $g(z)$ and for no one smaller value of $p$ we have $g(z+p)=g(z)$).

Thus certain analogy with the harmonic oscillator becomes clear. How can we use this same analogy bearing in mind, as this was discussed above, "Pólya – like" approach? Having no better ideas, in integral (2) we change $n+1/2$ to $\hat{H}_{osc}$ and introduce a unity operator $\hat{I}$ writing formally an operator equality:

$$\frac{-1-i}{2\sqrt{2\pi}} \int_{-1}^{\infty} \frac{\varsigma'}{\varsigma}(-2\hat{H}_{osc} + 5/4 + it)e^{-\pi t}dt +$$
$$\frac{(-1-i)e^{\pi}}{2\sqrt{2\pi}} \int_{0}^{\infty} \frac{\varsigma'}{\varsigma}(y - 2\hat{H}_{osc} + 5/4)(i\cos(\pi y) - \sin(\pi y))dy - (\Sigma_{\pi} - 3/2 + i)\hat{I} = (n+1/2)\hat{I} \quad (3).$$

Here integrands are to be understood as Taylor expansion: for any $\kappa$, where $\kappa = 5/4 + it$ or $\kappa = 5/4 + y$ with real $t \geq -1$, $y \geq 0$, we have $\frac{\varsigma'}{\varsigma}(\kappa - 2\delta) = \sum_{n=0}^{\infty} C_n(\kappa)\delta^n$, so write $\frac{\varsigma'}{\varsigma}(\kappa - 2\hat{H}_{osc}) = \sum_{n=0}^{\infty} C_n(\kappa)(\hat{H}_{osc}^n)$ quite similarly to e.g. the treatment of "an operator" $\exp(-a\hat{H})$ in many textbooks. Alternatively, we may introduce the functions

$$F(\delta) = \int_{-1}^{\infty} \frac{\varsigma'}{\varsigma}(5/4 + it - 2\delta)e^{-\pi t}dt, \qquad G(\delta) = \int_{0}^{\infty} \frac{\varsigma'}{\varsigma}(y + 5/4 - 2\delta)(i\cos(\pi y) - \sin(\pi y))dy \qquad \text{and}$$

consider their Taylor expansions instead. Operators should be used in some appropriate Hilbert functional space which we do not specify, implicitly supposing an analogy with the known case of the harmonic oscillator operator, its eigenfunctions and "superoperators", such as the aforementioned $\exp(-a\hat{H}_{osc})$.

What we have written up to now was, in a sense, "unconditionally true". Now our claim is that if and only if the operator

$$\hat{H}_{Polya} = \frac{-1-i}{2\sqrt{2\pi}} \int_{-1}^{\infty} \frac{\varsigma'}{\varsigma}(-2\hat{H}_{osc} + 5/4 + it)e^{-\pi t}dt +$$
$$\frac{(-1-i)e^{\pi}}{2\sqrt{2\pi}} \int_{0}^{\infty} \frac{\varsigma'}{\varsigma}(y - 2\hat{H}_{osc} + 5/4)(i\cos(\pi y) - \sin(\pi y))dy - (\Sigma_{\pi,RH} - 3/2 + i)\hat{I} \quad (4)$$

coincides with the harmonic operator Hamiltonian $\hat{H}_{osc}$, RH holds true, and to establish (better to say, "to illustrate") the possibility to construct such an operator, was the main aim of the paper.

We would like to finish this section remarking that the use of an infinite line $b_n - i$, $\infty - i$ is not necessary here, similar operator can be obtained exploiting any rectangular contour with the vertices at $b_n - i$, $b_n + iX$, $d + iX$, $d - i$ for any real $d>1$.

## 3. Conclusions

Certainly, a number of questions immediately arise here. For example, while the expansion $\frac{\varsigma'}{\varsigma}(\kappa - 2\delta) = \sum_{n=0}^{\infty} C_n(\kappa) \delta^n$ can be used, with an appropriate convergence radius, at any point along the lines $5/4 + it$, $5/4 + y$ (note that our integration path not only does not pass any peculiarity of the function $\varsigma'/\varsigma$ but also fully lies to the right of the point $z=1$), the questions of convergence of related "operator expansion" as well as its practical computation are difficult and we are not able to answer to them. We repeat that we do not pretend on a rigorous treatment of the problem here, our aim is different and more of heuristic nature: to try to propose a kind of new setting for the Hilbert – Pólya approach, as this has been discussed in the Introductory section. In this sense, we believe, this result might be useful, and we hope that another, more interesting and rigorously constructed "Pólya-like" operators of this type will follow.

S. Beltraminelli, CERFIM, Research Center for Mathematics and Physics, PO Box 1132, 6600 Locarno, Switzerland.
E-mail: **Stefano.beltraminelli@ti.ch**

D. Merlini, CERFIM, Research Center for Mathematics and Physics, PO Box 1132, 6600 Locarno, Switzerland.
E-mail: **merlini@cerfim.ch**

S. K. Sekatskii, Laboratoire de Physique de la Matière Vivante, IPSB, BSP 408, Ecole Polytechnique Fédérale de Lausanne, CH1015 Lausanne-Dorigny, Switzerland.
E-mail : **serguei.sekatski@epfl.ch**